\documentclass[11pt]{article}
\usepackage{amsmath}
\usepackage{amssymb}

\begin{document}
\renewcommand{\thesection}{\arabic{section}}
\renewcommand{\theequation}{\thesection.\arabic{equation}}
\title{{\Large   Existence of Positive Solutions for a class of Quasilinear Schr\"{o}dinger Equations of Choquard type}
\thanks{This work was supported partially by the National
Natural Science Foundation of China (11562021). }}
\author{{Shaoxiong Chen}
\thanks{gxmail@126.\ com} \ and \  {Xian Wu}
\thanks{Corresponding Author:\ wuxian2042@163.\ com}\\
{\small Department of Mathematics,\ Yunnan Normal University,\
Kunming, Yunnan 650092, P. R. China}}
\date{}
\maketitle \textbf{Abstract:} \ \ In this paper, we study the
following quasilinear Schr\"{o}dinger equation of Choquard type
$$
-\triangle u+V(x)u-\triangle
(u^{2})u=(I_\alpha *|u|^p)|u|^{p-2}u, \ \  x \in
\mathbb{R}^{N},
$$
 where $N\geq 3$,\ $0<\alpha<N$,
 $\frac{2(N+\alpha)}{N}\leq p<\frac{2(N+\alpha)}{N-2}$ and $I_\alpha$ is a Riesz potential. Under appropriate assumptions on
$V(x)$, we establish the existence of positive solutions.
\par

\textbf{Key Words}: quasilinear Schr\"{o}dinger equations,
Choquard type, Variational methods.\par

\section*{1. Introduction and Preliminaries}
\ \ \ \ Consider the following quasilinear Schr\"{o}dinger equation
of Choquard type
$$
-\triangle u+V(x)u-\triangle
(u^{2})u=(I_\alpha *|u|^p)|u|^{p-2}u, \ \  x \in
\mathbb{R}^{N},
 \ \eqno(1. 1) $$ where $N\geq 3$,\ $0<\alpha<N$,
 $\frac{2(N+\alpha)}{N}\leq p<\frac{2(N+\alpha)}{N-2}$, $V\in C(\mathbb{R}^{N},
\mathbb{R})$ and $I_\alpha:\ \mathbb{R}^{N}\rightarrow \mathbb{R}$ is the Riesz potential defined by
$$
I_\alpha (x)=\frac{\Gamma\left(\frac{N-\alpha}{2}\right)}{\Gamma\left(\frac{\alpha}{2}\right)\pi^{N/2}2^\alpha|x|^{N-\alpha}}
:=\frac{A_\alpha}{|x|^{N-\alpha}},
$$
and $\Gamma$ is the Gamma function.

 Eq.(1.1) is related to quasilinear Schr\"{o}dinger
equations of more general form
$$i\psi_t+\triangle\psi-V(x)\psi+k\triangle(h(|\psi|^2))h'(|\psi|^2)\psi+g(x,\psi)=0 \ \eqno(1.2)$$
where $V=V(x)$, $x\in R^N$, is a given potential, $k$ is a real constant and $h$,
$g$ are real functions. The quasilinear schr\"{o}dinger equations
(1.2) are derived as models of several physical phenomena, such as
see \cite{SK, EKL, AA, AN, MM}. It has received considerable attention in mathematical analysis during the last 10 years \cite{MP4}. When $k\neq 0$ and $g$ is a local term, several methods can be used to solve Eq. $(1.2)$. The existence of a positive
ground state solution has been proved in \cite{JZ, MKZ} by using a
constrained minimization approach. The problem is transformed to a
semilinear one  in \cite{ML, JYZ2, DMS} by a change of variables. Nehari method is used to get the existence results of
ground state solutions in \cite{JYZ, DG}.\par

For $k=0$, $V(x)\equiv 1$ and $g(x,\psi)=\left(I_\alpha *|\psi|^p\right)|\psi|^{p-2}\psi$, Eq.$(1.2)$
 is usually called the nonlinear Choquard or Choquard-Pekar
equation. It has several physical origins. The problem appeared at
least as early as in 1954, in a work by S. I. Pekar describing
the quantum mechanics of a polaron at rest. In 1976, for $N=3, \alpha=2$, P. Choquard used it to describe
an electron trapped in its own hole, in a certain approximation to Hartree-Fock theory of
one component plasma(see \cite{VJ}). In that case, a ground state solution is found in \cite{VJ}.
For $p\in \left(\frac{N+\alpha}{N},\frac{N+\alpha}{N-2}\right)$,
Moroz and Van Schaftingen in \cite{VJ} proved existence, qualitative properties and decay estimates of ground
state solutions. In \cite{J},
Seok consider a critical version of nonlinear Choquard equation
$$
-\triangle u+u=(I_\alpha *|u|^p)|u|^{p-2}u+\lambda |u|^{2^*-2}u, \ \  x \in
\mathbb{R}^{N}.
$$
Using some perturbation arguments, Seok gets a family of nontrivial solutions. It converges to a least energy solution
of the limiting critical local problem as $\alpha\rightarrow 0$.
We can see some works of literature about the above equation in $\cite {CAM, SJX, HY}$.

 In this paper, we study
the existence of positive solution of $(1.1)$ by  a variational
argument. We need the following several notations. If $x\in \mathbb{R}^N$ and $R>0$, the closed ball with
center at $x$ and radius $R$ is always denoted by $B_R(x)$. Let
$C_{0}^{\infty}(\mathbb{R}^{N})$ be the collection of smooth
functions with compact support. For $N\geq 3$, let
$$D^{1, 2}(\mathbb{R}^{N}):=\{u\in L^{2^*}(\mathbb{R}^{N}):\nabla
u\in L^{2}(\mathbb{R}^{N}) \}$$ with the norm
$$\|u\|_{D^{1, 2}}^2=\int_{\mathbb{R}^{N}}|\nabla u|^2.$$
By the Sobolev inequality, $D^{1, 2}(\mathbb{R}^{N})$ is continuously
embedded in $L^{2^*}(\mathbb{R}^{N})$. Let
$$H^{1}(\mathbb{R}^{N}):= \{u \in L^{2}(\mathbb{R}^{N}): \nabla u
\in L^{2}(\mathbb{R}^{N}) \}$$ with the inner product
$$\langle u, v\rangle_{H^1}=\int_{\mathbb{R}^{N}}(\nabla u\cdot\nabla
v+uv)$$
 and the norm
 $$\|u\|:=\|u\|_{H^1}=\langle u, u\rangle_{H^1}^{1/2}.$$
We denote the norm of $L^q(\mathbb{R}^{N})$ by $|\cdot|_q$.

In the following, we always assume $V\in C(\mathbb{R}^{N},
\mathbb{R})$ and $\inf\limits_{\mathbb {R}^N}V(x)\geq V_0>0$. Let us
consider the following two assumptions:\par

 $(V_1)$\ \ $V(x)$ is periodic in
each variable of $x_1,\ \cdots,\ x_N$.\par
$(V_2)$\ \ $V(x)\leq V_{\infty}:=\lim\limits_{|y|\rightarrow
 \infty}V(y)<\infty$ for all $x\in \mathbb{R}^{N}$ and $V_0<V_\infty$.\par

Eq. $(1. 1)$ is the Euler-Lagrange equation  of the energy
functional
$$J(u)  =
\frac{1}{2}\int_{\mathbb{R}^{N}}[(1+2u^{2})|\nabla u|^{2} +
V(x)u^{2}] -\frac{1}{2p}
\int_{\mathbb{R}^{N}}\left(I_\alpha *|u^+|^p\right)|u^+|^p,$$
where $u^+=\max\{u,0\}$.

The main result of this paper is stated as follows: \par

 {\bf Theorem 1.1.} \ \ {\it Suppose that $N\geq 3$,
$\frac{2(N+\alpha)}{N}\leq p<\frac{2(N+\alpha)}{N-2}$ and the potential function $V$ satisfies condition $(V_1)$ or $(V_2)$. Then Eq. (1.1) possesses a positive solution $u\in H^1(\mathbb{R}^N)$.} \par

 {\bf Remark 1.2.} \ \ By the Hardy-Littlewood-Sobolev inequality and the Sobolev embedding theorem,  the natural interval for considering the Choquard equation is $[\frac{2(N+\alpha)}{N},\frac{2(N+\alpha)}{N-2}]$, however, the critical case $p=\frac{2(N+\alpha)}{N-2}$ is not considered in Theorem 1.1.\par
   \ \ To deal with this type of problem, difficulties lie in two aspects. First, the approach of proving Theorem 1.1 is inspired by \cite{ML,JOS,XA,JYZ2}. Because the nonlinearity of Eq. $(1.1)$ is nonlocal, the techniques in these papers developed for the local case cannot be adopted directly. Second, another major difficulty here is that the energy functional $J(u)$ is not well
defined for all $u\in H^1(\mathbb{R}^{N})$ if $N\geq 3$. To overcome these difficulties, we need to make a change $u=f(v)$ of variable used in
\cite{JYZ2}, to analyze other properties of $f$ more deeply due to the nonlocal nonlinearity and to develop some different techniques. \par
  Following the idea in \cite{JYZ2}, let $f$ be defined by
$$f'(t)=\frac{1}{\sqrt{1+2f^{2}(t)}}$$ on $[0, +\infty)$, $f(0)=0$
and $f(-t)=-f(t)$ on $(-\infty, 0]$.  Then $f$ has following
properties (see \cite{YX}):\par

 $(f_1)$ \ $f$ is
uniquely defined $C^{\infty}$ function and invertible.\par

$(f_2)$ \ $0<f'(t)\leq 1$  for all $t\in \mathbb{R}$.\par

$(f_3)$ \ $|f(t)|\leq |t|$ for all $t\in \mathbb{R}$.\par

$(f_4)$ \  $\frac{1}{2}f(t)\leq tf^{\prime}(t)\leq f(t)$ \ for all
$t\geq 0$ and $f(t)\leq tf^{\prime}(t)\leq\frac{1}{2}f(t)$ for all
$t\leq 0$.\par

$(f_5)$ \ $|f(t)|\leq 2^{\frac{1}{4}}|t|^{\frac{1}{2}}$  for all $t\in
\mathbb{R}$.\par

$(f_6)$ \ There exists a positive constant $C$ such that
$$\begin{array}{ll}
|f(t)|\geq  \left\{
\begin{array}{ll}
C|t|, &   |t|\leq 1,\\
C|t|^{\frac{1}{2}}, & |t|\geq 1.
\end{array}
\right.
\end{array}$$
\par

$(f_{7})$ \ $|f(t)f^{\prime}(t)|\leq \frac{1}{\sqrt{2}}$ for all
$t\in \mathbb{R}$.\par

$(f_{8})$ \ For each $\xi>0$, there exists $C(\xi)>0$ such that
$f^2(\xi t)\leq C(\xi)f^2(t)$.\par

After the change $u=f(v)$ of variable,  Eq. $(1. 1)$ can be rewritten as
$$
-\triangle v+V(x)f(v)f'(v)=(I_\alpha *|f(v)|^p)|f(v)|^{p-2}f(v)f'(v), \ \  x \in
\mathbb{R}^{N},
 \ \eqno(1.3)
$$
and $J(u)$ can be reduced to
$$I(v):=\frac{1}{2}\int_{\mathbb{R}^{N}}(|\nabla v|^{2} +
V(x)f^{2}(v))-
\frac{1}{2p}
\int_{\mathbb{R}^{N}}\left(I_\alpha *|f(v^+)|^p\right)|f(v^+)|^p.
\eqno(1.4)$$

From $(V_1)$ (or $(V_2)$ ), $(f_3)$, $(f_7)$ and the Hardy-Littlewood-Sobolev inequality (Lemma 2.1),
we can deduce that the functional $I\in C^1(H^1(\mathbb{R}^N))$ for $\frac{2(N+\alpha)}{N}\leq p\leq\frac{2(N+\alpha)}{N-2}$.\par

It is easy to see that if  $v\in H^1(\mathbb{R}^N)$ is a critical point of
 $I$, i.e.,
$$
\begin{aligned}
\langle I'(v), \varphi\rangle=\int_{\mathbb{R}^{N}}\nabla v\nabla
\varphi + \int_{\mathbb{R}^{N}}V(x)f(v)f'(v)\varphi
-\int_{\mathbb{R}^{N}}\left(I_\alpha *|f(v^+)|^p\right)|f(v^+)|^{p-1}f'(v^+)\varphi
,
\end{aligned}
$$
 for all $\varphi \in C^\infty_0(\mathbb{R}^N)$, then
 $v$ is a weak solution of Eq. $(1. 3)$ and
$u:=f(v)$ is a weak solution of Eq. (1.1).

We use $C$ or $C_i$ to denote various positive
 constants in context. The outline of the paper is as follows. In Section 2, we give some preliminary results and the regularity of solutions. In Section 3, we prove Theorem 1.1 by using the mountain pass theorem.\par

\section*{2. Some lemmas}

{\bf Lemma 2.1} (\cite{VJ}) \ (Hardy-Littlewood-Sobolev inequality).\ \ {\it Let $r, s>1$ and $0<\alpha <N$
be such that
$$
\frac{1}{r}+\frac{1}{s}-\frac{\alpha}{N}=1.
$$
Let $g\in L^r(\mathbb{R}^N)$ and $h\in L^s(\mathbb{R}^N)$. There exists a sharp constant $C(r, s, N, \alpha)$,
independent of $g$,\ $h$,\ such that
$$
\int_{\mathbb{R}^{N}}\int_{\mathbb{R}^{N}}\frac{g(x)h(y)}{|x-y|^{N-\alpha}}\leq C(r, s, N, \alpha)|g|_r|h|_s.
 \ \eqno(2.1)
$$
}
{\bf Remark 2.2} (1)\ Hardy-Littlewood-Sobolev inequality can be also stated that for every
$s\in (1, \frac{N}{\alpha})$, for every $v\in L^s(\mathbb{R}^{N})$, $I_\alpha *v\in L^{\frac{Ns}{N-\alpha s}}(\mathbb{R}^{N})$ and the following inequality holds:
$$
\int_{\mathbb{R}^{N}}\int_{\mathbb{R}^{N}}|I_\alpha *v|^{\frac{Ns}{N-\alpha s}}\leq
C\left(\int_{\mathbb{R}^{N}}|v|^s\right)^{\frac{N}{N-\alpha s}},
 \ \eqno(2.2)
$$
where $C>0$ depends only on $\alpha, N$ and $s$.\par

(2) By Lemma 2.1 and $(f_5)$,
$$
\int_{\mathbb{R}^{N}}\int_{\mathbb{R}^{N}}\frac{|f(v(x))|^p|f(v(y))|^p}{|x-y|^{N-\alpha}}\leq
C_1\int_{\mathbb{R}^{N}}\int_{\mathbb{R}^{N}}\frac{|v(x)|^{\frac{p}{2}}|v(y)|^{\frac{p}{2}}}{|x-y|^{N-\alpha}}\leq C|v|_{\frac{p}{2}r}^p
$$
if $v\in L^{\frac{p}{2}r}(\mathbb{R}^N)$ for $r>1$ with
$$
\frac{2}{r}-\frac{\alpha}{N}=1.
$$
Since we will work on $H^1(\mathbb{R}^N)$, by the Sobolev embedding theorem,\ we must require that $\frac{p}{2}r\in [2, 2^*]$. We have
$$
\frac{2(N+\alpha)}{N}\leq p\leq\frac{2(N+\alpha)}{N-2}.
$$
\ \ \ \ (3) $H^1(\mathbb{R}^N)\subset L^{\frac{2Nq}{N+\alpha}}(\mathbb{R}^N)$ if and only if $\frac{N+\alpha}{N}\leq q\leq\frac{N+\alpha}{N-2}$ where $q:=\frac{p}{2}$ (see $\cite {VJ}$).

{\bf Lemma 2.3}\ \ {\it $|f(t)|\leq 2^{\frac{1}{4}}|t|^l$ for all $l\in [\frac{1}{2}, 1]$, $t\in \mathbb{R}$,.}

{\bf Proof.} By $(f_3), (f_5)$,\ for $|t|\leq 1$,
 $$|f(t)|\leq |t|\leq |t|^l\leq 2^{\frac{1}{4}}|t|^l.$$
For $|t|>1$,
$$|f(t)|\leq 2^{\frac{1}{4}}|t|^{\frac{1}{2}}\leq 2^{\frac{1}{4}}|t|^l.\ \ \Box$$

{\bf Lemma 2.4}\ \ {\it For $N\geq 3,\ 0<\alpha<N$,\ $\frac{2(N+\alpha)}{N}\leq p\leq\frac{2(N+\alpha)}{N-2}$, there exists $C>0$ and $r\in [2, 2^*)$ such that
$$
|f(t)|^p\leq Ct^r,\ |f(t)|^{p-2}f(t)f'(t)\leq Ct^{r-1},\ \ \forall \ t\geq 0.
 \ \eqno(2.3)
$$
}

{\bf Proof.}  If $\frac{p}{2}\geq 2$, we put $r=\frac{p}{2}$. Since $p\leq\frac{2(N+\alpha)}{N-2}$ and $\alpha<N$, $r\in [2, 2^*)$.  By $(f_5),\ (f_{7})$, we have
$$
|f(t)|^{p-2}f(t)f'(t)\leq Ct^{\frac{1}{2}(p-2)}=Ct^{r-1}
$$
and
$$|f(t)|^p=|f^2(t)|^{\frac{p}{2}}\leq Ct^r$$
for all $t\geq 0$.\par

If $\frac{p}{2}<2$,\ by $2<\frac{2(N+\alpha)}{N}\leq p$ and Lemma 2.3, then there exists $l\in [\frac{1}{2},1]$ such that $r:=lp=2$.
Hence
$$|f(t)|^p\leq C|t|^{lp}=Ct^{r}, \ \ \ \ \forall \ t\geq 0.$$
Moreover,
for $t\geq 1$, using $(f_{7})$,
$$
|f(t)|^{p-2}f(t)f'(t)\leq Ct^{l(p-2)}\leq Ct^{l(p-2)}\leq Ct^{lp-1}=Ct^{r-1}.
$$
For $0\leq t<1$, by $(f_2)$ and $(f_3)$, we have
$$
|f(t)|^{p-2}f(t)f'(t)\leq Ct^{l(p-2)}t\leq Ct^{lp-2}t=Ct^{r-1}.\ \ \Box
$$

In the spirit of the argument developed by the Proposition 4.1 in \cite {VJ} and the Lemma 2.1 in \cite {CM},\ we have the following Lemma 2.5.

{\bf Lemma 2.5}\ \ {\it Assume that either assumption $(V_1)$ or $(V_2)$ holds. Let $N\geq 3,\ 0<\alpha<N,\ \frac{2(N+\alpha)}{N}\leq p<\frac{2(N+\alpha)}{N-2}$,
$v\in H^1(\mathbb{R}^N)$ be a nontrivial and nonnegative weak solution of $(1.3)$, then the following properties hold:\par

(1)\ For every $s\in \mathbb{R}$ with
$\frac{\alpha}{N}\left(1-\frac{2}{p}\right)-\frac{2}{N}<\frac{1}{s}<1$, we have
$v\in L^s_{loc}(\mathbb{R}^N)$.\par

(2)\ $I_\alpha *|f(v)|^p\in L^\infty_{loc} (\mathbb{R}^N)$.\par

(3)\ $v\in W^{2,\ q}_{loc}(\mathbb{R}^N)$ for every $q\in [2, \infty)$.\par

(4)\ For any $\lambda\in (0, 1)$ such that
$$
\left(I_\alpha *|f(v)|^p\right)|f(v)|^{p-2}f(v)f'(v)\in C^{0, \lambda}_{loc}(\mathbb{R}^N).
$$

(5)\ $v$ is of class $C^{1, \lambda}$ for every $\lambda\in (0, 1)$ and $v>0$.}

{\bf Proof.\ (1)} We use an iterating argument (see $\cite {VJ}$). Put $r:=\frac{p}{2}>1 $, since $v\in H^1(\mathbb{R}^N)$, By Remark 2.2-(3), $v\in L^{\frac{2Nr}{N+\alpha}}(\mathbb{R}^N)$. Set $s_0=\underline{s}_0=\overline{s}_0:=\frac{2Nr}{N+\alpha}$ and
$$
\frac{1}{t}=\frac{r}{s_0}-\frac{\alpha}{N}\in (0, 1).
$$
We have $I_\alpha *|v|^r\in L^t(\mathbb{R}^N)$. Further, set
$$
\frac{1}{q_{_0}}=\frac{2r-1}{s_0}-\frac{\alpha}{N}=\frac{r-1}{s_0}+\frac{1}{t}\in (0, 1).
$$
We get $\left(I_\alpha *|v|^r\right)|v|^{r-1}\in L^{q_{_0}}(\mathbb{R}^N)$. By $(f_4), \ (f_5), \ (f_{7}), \ v\geq 0$ and
$$
\left(I_\alpha *|f(v)|^p\right)|f(v)|^{p-2}f(v)f'(v)
\leq C\left(I_\alpha *|v|^{\frac{p}{2}}\right)|v|^{\frac{p-2}{2}}= C\left(I_\alpha *|v|^r\right)|v|^{r-1},
$$
one has $$\left(I_\alpha *|f(v)|^p\right)|f(v)|^{p-2}f(v)f'(v)\in L^{q_{_0}}(\mathbb{R}^N).$$
Set
$$
\begin{array}{ll}
c(x)=\left\{
\begin{array}{ll}
\frac{V(x)f(v(x))f'(v(x))}{v(x)}, &   v(x)\neq 0,\\
1, & v(x)=0.
\end{array}
\right.
\end{array}
$$
Using $(f_2), \ (f_3), \ (V_1)$ (or $(V_2)$), we know that $c(x)\in L^\infty_{loc}(\mathbb{R}^N)$ and $c(x)v=V(x)f(v)f'(v)$.
Hence $v$ is a weak solution of the following equation
$$
-\triangle v+c(x)v=(I_\alpha *|f(v)|^p)|f(v)|^{p-2}f(v)f'(v), \ \  x \in
\mathbb{R}^{N},
$$
and hence, $v\in W^{2, {q_{_0}}}_{loc}(\mathbb{R}^N)$ by the Theorem 9.1.4 in \cite {WYW}. \par

Notice that $\frac{2(N+\alpha)}{N}\leq p<\frac{2(N+\alpha)}{N-2}$.
By the Sobolev embedding theorem,  $v\in L^s_{loc}(\mathbb{R}^N)$
provided
$$
\frac{2r-1}{s_0}-\frac{\alpha+2}{N}\leq\frac{1}{s}
\leq\frac{2r-1}{s_0}-\frac{\alpha}{N}.
$$
Furthermore,
$$
\frac{2r-1}{\overline{s}_0}-\frac{\alpha+2}{N}<\frac{1}{\overline{s}_0}=\frac{1}{s_0}=\frac{1}{\underline{s}_0}
\leq\frac{2r-1}{\underline{s}_0}-\frac{\alpha}{N}.
$$
Set $\frac{1}{\underline{s}_1}=\frac{2r-1}{\underline{s}_0}-\frac{\alpha}{N}$ and  $\frac{1}{\overline{s}_1}=\frac{2r-1}{\overline{s}_0}-\frac{\alpha+2}{N}$.
For every $s\in (\underline{s}_1, \overline{s}_1)$, we have $v\in L^s_{loc}(\mathbb{R}^N)$. Then we can prove
$v\in W^{2, {q}}_{loc}(\mathbb{R}^N)$
 provided
$$
\frac{2r-1}{\overline{s}_1}-\frac{\alpha}{N}<\frac{1}{q}
<\frac{2r-1}{\underline{s}_1}-\frac{\alpha}{N}
$$
and
$$
\frac{\alpha}{N}\left(1-\frac{1}{r}\right)<\frac{1}{q}<1.
$$
In fact, for every $s$ with $\frac{1}{q}=\frac{2r-1}{s}-\frac{\alpha}{N}$, the above inequalities follows that
$$
\frac{\alpha}{Nr}<\frac{1}{s}<\frac{N+\alpha}{N(2r-1)}
$$
and
$$s\in (\underline{s}_1, \overline{s}_1).$$
Hence
$$
\frac{1}{t}=\frac{r}{s}-\frac{\alpha}{N}\in (0, 1),
$$
and hence $I_\alpha *|v|^r\in L^t_{loc}(\mathbb{R}^N)$. Further, since
$$
\frac{1}{q}=\frac{2r-1}{s}-\frac{\alpha}{N}=\frac{r-1}{s}+\frac{1}{t}\in (0,1),
$$
one has $\left(I_\alpha *|v|^r\right)|v|^{r-1}\in L^{q}_{loc}(\mathbb{R}^N)$. Similar to the above segment, we can prove $v\in W^{2, {q}}_{loc}(\mathbb{R}^N)$. \par

Set
$$
\left(\frac{1}{\overline{q}_1}, \frac{1}{\underline{q}_1}\right)=\left(\frac{2r-1}{\overline{s}_1}-\frac{\alpha}{N},
\frac{2r-1}{\underline{s}_1}-\frac{\alpha}{N}\right)\cap \left(\frac{\alpha}{N}\left(1-\frac{1}{r}\right), 1\right) (\neq \emptyset),
$$
If $\frac{1}{\overline{q}_1}=\frac{\alpha}{N}\left(1-\frac{1}{r}\right)$ and $\frac{1}{\underline{q}_1}=\frac{2r-1}{\underline{s}_1}-\frac{\alpha}{N}<1$,  we use $\underline{q}_1$
in replace of $q_0$ in the above argument. Similarly, if $\frac{1}{\overline{q}_1}=\frac{2r-1}{\overline{s}_1}-\frac{\alpha}{N}>\frac{\alpha}{N}\left(1-\frac{1}{r}\right)$ and $\frac{1}{\underline{q}_1}=1$,  we use $\overline{q}_1$ in replace of $q_0$.
So, without loss of generality, we can assume that $\frac{\alpha}{N}\left(1-\frac{1}{r}\right)<\frac{2r-1}{\overline{s}_1}-\frac{\alpha}{N}$ and
$\frac{2r-1}{\underline{s}_1}-\frac{\alpha}{N}<1$ are satisfied. Hence, for every $q\in \mathbb{R}$ with $\frac{2r-1}{\overline{s}_1}-\frac{\alpha}{N}<\frac{1}{q}
<\frac{2r-1}{\underline{s}_1}-\frac{\alpha}{N}$, we know that  $v\in  W^{2, {q}}_{loc}(\mathbb{R}^N)$. \par

By the Sobolev embedding theorem again , $v\in L^s_{loc}(\mathbb{R}^N)$
provided
$$
\frac{2r-1}{\overline{s}_1}-\frac{\alpha+2}{N}\leq\frac{1}{s}
\leq\frac{2r-1}{\underline{s}_1}-\frac{\alpha}{N}.
$$
We also have
$$
\frac{\alpha}{N}\left(1-\frac{2}{p}\right)-\frac{2}{N}<\frac{1}{s}<1.
$$
Furthermore, by $\frac{2(N+\alpha)}{N}\leq p<\frac{2(N+\alpha)}{N-2}$,
$$
\frac{2r-1}{\overline{s}_1}-\frac{\alpha+2}{N}<\frac{1}{\overline{s}_1}<\frac{1}{\overline{s}_0}=\frac{1}{s_0}=\frac{1}{\underline{s}_0}\leq\frac{1}{\underline{s}_1}
\leq\frac{2r-1}{\underline{s}_1}-\frac{\alpha}{N}.
$$
Set $\frac{1}{\underline{s}_2}=\frac{2r-1}{\underline{s}_1}-\frac{\alpha}{N}$ and  $\frac{1}{\overline{s}_2}=\frac{2r-1}{\overline{s}_1}-\frac{\alpha+2}{N}$. Then, for each $s\in (\underline{s}_2, \overline{s}_2)$, we have $v\in L^s_{loc}(\mathbb{R}^N)$. Continue the above process by setting
$$
\frac{1}{\overline{s}_{k+1}}=\frac{2r-1}{\overline{s}_k}-\frac{\alpha+2}{N},\ \
\frac{1}{\underline{s}_{k+1}}=\frac{2r-1}{\underline{s}_k}-\frac{\alpha}{N}.
$$
Notice that
$$
\frac{1}{\overline{s}_{k+1}}-\frac{1}{\overline{s}_k}=(2r-1)^k\left(\frac{1}{\overline{s}_1}-\frac{1}{\overline{s}_0}\right)<0,\ \
\frac{1}{\underline{s}_{k+1}}-\frac{1}{\underline{s}_k}=(2r-1)^k\left(\frac{1}{\underline{s}_1}-\frac{1}{\underline{s}_0}\right)
\geq 0.
$$
If $r=p/2=\frac{N}{N+\alpha}$, then $\frac{1}{\underline{s}_{k+1}}=\frac{1}{\underline{s}_k}=\cdots=\frac{1}{\underline{s}_0}=\frac{1}{s}<1$.
If $\frac{2(N+\alpha)}{N}<p<\frac{2(N+\alpha)}{N-2}$, then $\frac{1}{\underline{s}_1}-\frac{1}{\underline{s}_0}>0$.
Since $r>1$, there exists $k_0\geq 1$ such that
$$
\frac{1}{\overline{s}_k}<\frac{\alpha}{N}\left(1-\frac{2}{p}\right)-\frac{2}{N}<1<\frac{1}{\underline{s}_k}.
$$
Hence (1) holds.\par

{\bf (2)} Since $p<\frac{2(N+\alpha)}{N-2}$,
$$
\frac{\alpha}{N}\left(1-\frac{2}{p}\right)-\frac{2}{N}<\frac{\alpha^2-2N}{N(N+\alpha)}<\frac{\alpha^2-2\alpha}{N(N+\alpha)}<\frac{2\alpha}{Np}.
$$
Fix $\frac{1}{q}\in \left(\frac{\alpha}{N}\left(1-\frac{2}{p}\right)-\frac{2}{N}, \frac{2\alpha}{Np}\right)$. It implies that $v\in L^q_{loc}(\mathbb {R}^N)$ and $N-1+\frac{2q}{2q-p}(\alpha-N)>-1$. By virtue of the Lemma 2.4 and the H\"older's inequality, for every closed ball $B\subset \mathbb{R}^N$, there exists $r\in [2, 2^*)$ such that
$$\begin{aligned}
\left(I_\alpha *|f(v)|^p\right)(x)
& =A_\alpha\left(\int_{|x-y|>1}\frac{|f(v(y))|^p}{|x-y|^{N-\alpha}}dy+\int_{|x-y|\leq 1}\frac{|f(v(y))|^p}{|x-y|^{N-\alpha}}dy\right)\\
& \leq
C_1\left(\int_{|x-y|>1}\frac{|v(y)|^r}{|x-y|^{N-\alpha}}dy+\int_{|x-y|\leq 1}\frac{|f(v(y))|^p}{|x-y|^{N-\alpha}}dy\right)\\
& \leq
C_1\int_{|x-y|>1}|v(y)|^rdy+C_3\int_{|x-y|\leq 1}\frac{|v(y)|^{\frac{p}{2}}}{|x-y|^{N-\alpha}}dy\\
& \leq
C_2+C_3\left(\int_{|x-y|\leq 1}|x-y|^{(\alpha-N)\cdot\frac{2q}{2q-p}}dy\right)^{1-\frac{p}{2q}}\left(\int_{|x-y|\leq 1}|v(y)|^{\frac{p}{2}\cdot\frac{2q}{p}}dy\right)^{\frac{p}{2q}}\\
& \leq C
\end{aligned}$$
for all $x\in B$.
That is  $I_\alpha *|f(v)|^p\in L^\infty_{loc} (\mathbb{R}^N)$.\par

{\bf (3)} For every closed ball $B\subset \mathbb{R}^N$, set
$$
a(x)=\left|\frac{\left(I_\alpha *|f(v)|^p\right)|f(v)|^{p-2}f(v)f'(v)-V(x)f(v)f'(v)}{1+|v|}\right|,\ \forall x\in B.
$$
By the Lemma 2.4, the conclusion (2), $(V_1)$ or $(V_2)$, it follows that
$$
a(x)\leq V(x)+C_1\frac{|f(v)|^{p-2}f(v)f'(v)}{1+|v|}\leq C_2+C_1|v|^{r-2},\ \forall x\in B.
$$
Since
$$\begin{aligned}
\int_{B}|v|^{(r-2)\cdot\frac{N}{2}}dx& =
\int_{B\cap \{|v|<1|\}}|v|^{(r-2)\cdot\frac{N}{2}}dx+\int_{B\cap \{|v|\geq 1\}}|v|^{(r-2)\cdot\frac{N}{2}}dx\\
& \leq
C_1+\int_{B\cap \{|v|\geq 1\}}|v|^{(\frac{2N}{N-2}-2)\cdot\frac{N}{2}}dx\\
& \leq C
\end{aligned}$$
for all $x\in B$,  $a(x)\in L^{\frac{N}{2}}(B)$. Hence $a(x)\in L^{\frac{N}{2}}_{loc}(\mathbb{R}^N)$. Apply the Lemma B.3 in $\cite{M}$ to the following equation
$$
-\Delta v=a(x)v.
$$
It follows that $v\in L^q_{loc}(\mathbb{R}^N)$ for any $q\in [2, +\infty)$. Hence, by the Theorem 9.1.4 of \cite {WYW}, $v\in W^{2, q}_{loc}(\mathbb{R}^N)$ for any $q\in [2, +\infty)$.\par

{\bf (4)} From (3) and the Sobolev embedding theorem, we have that $v\in C^{1,\lambda}_{loc}(\mathbb{R}^N)$ for any $\lambda\in (0,1)$.
Fix $R>0$. Take $\beta\in C^\infty_0(\mathbb{R}^N)$ such that $\beta(t)=1$ for $t\leq R$,\ $\beta(t)=0$ for $t\geq 2R$.
Write $I_\alpha=\beta I_\alpha+(1-\beta)I_\alpha$. Put
$$\begin{array}{ll}
g(x)=\left\{
\begin{array}{ll}
0, &   |x|\leq R,\\
(1-\beta(x))I_\alpha(x), & |x|\geq R.
\end{array}
\right.
\end{array}$$
Then $g(x)\in C_0^\infty (\mathbb{R}^{N})$.
By the Lemma 2.4,\ there exists $r\in [2,2^*)$ such that
$$
g'(x-\cdot)|f(v(\cdot))|^p\leq C|v(\cdot)|^r\in L^1(\mathbb{R}^{N})
$$
for every $x\in \mathbb{R}^{N}$.
Using the Lebesgue's dominated convergence theorem, we have
$$
\left(((1-\beta)I_\alpha)*|f(v)|^p\right)(x)=\int_{\mathbb{R}^{N}}g(x-y)|f(v(y))|^pdy\in C^\infty (\mathbb{R}^{N}).
$$
It follows from $p>2$, $f\in C^\infty(\mathbb{R})$ that $|f(v)|^p, |f(v)|^{p-2}f(v)f'(v)\in C^{0, \lambda}_{loc}(\mathbb{R}^N)$.
It is easy to see that $\beta I_\alpha\in L^1(\mathbb{R}^{N})$. Hence
$$
\left(\beta I_\alpha\right)*|f(v)|^p\in C^{0, \lambda}_{loc}(\mathbb{R}^N).
$$
Therefore,
$$
\left(I_\alpha *|f(v)|^p\right)|f(v)|^{p-2}f(v)f'(v)\in C^{0, \lambda}_{loc}(\mathbb{R}^N).
$$

{\bf (5)} It follows from (4) that $v\in C_{loc}^{1, \lambda}(\mathbb{R}^{N})$ for every $\lambda\in (0, 1)$. In the proof of conclusion (1), Eq. $(1.3)$ can be rewritten as
$$
-\triangle v+c(x)v=(I_\alpha *|f(v)|^p)|f(v)|^{p-2}f(v)f'(v).
$$
Since $v$ is a continuous nonnegative function, there exists $C>0$ such that $|c(x)|\leq C$. Using  the Theorem 8.19 in \cite{DN}, we get that $v>0$ in $\mathbb{R}^N$.\
$\Box$

\section*{3. Proof of Theorem 1.1}

 \ \ \ \ We prove that the functional $I$ exhibits the mountain pass geometry.

{\bf Lemma 3.1.}\ \ {\it There exist $\rho_{_{_0}}, \alpha>0$ such that
$$
I(v)\geq \alpha, \ \ \ \forall \  v \in \{v\in H^1(\mathbb{R}^N): \|v\|=\rho_{_{_0}}\}.
$$
}
{\bf Proof.} From \cite{XA} we get that
there exist $C_1>0, \rho_{_{_1}}>0$ such that
$$
\int_{\mathbb{R}^N}\left(|\nabla v|^2+V(x)f^2(v)\right)\geq C_1\|v\|^2
 \ \eqno(3.1)
$$
whenever $\|v\|\leq \rho_{_{_1}}$. The above inequality was derived in \cite{XA} for $(V_1)$. Checking the proof of \cite{XA}, we know that this inequality holds for $(V_2)$, too. Notice that $\frac{Np}{N+\alpha}\in [2, 2^*)$. By $(f_5)$, (2.1), (3.1) and
the Sobolev embedding theorem, we get
$$\begin{aligned}
I(v)& \geq
\frac{C_1}{2}\|v\|^2-\frac{1}{2p}\int_{\mathbb{R}^N}\left(I_\alpha*|f(v^+)|^p\right)|f(v^+)|^p\\
& \geq
\frac{C_1}{2}\|v\|^2-C_2\int_{\mathbb{R}^N}\left(I_\alpha*|v|^{p/2}\right)|v|^{p/2}\\
& \geq
\frac{C_1}{2}\|v\|^2-C_2\left(\int_{\mathbb{R}^N}|v|^{\frac{Np}{N+\alpha}}\right)^{\frac{N+\alpha}{N}}\\
& \geq
\frac{C_1}{2}\|v\|^2-C_3\|v\|^p\\
& \geq
\|v\|^2\left(\frac{C_1}{2}-C_3\|v\|^{p-2}\right)
\end{aligned}$$
whenever $\|v\|\leq \rho_{_{_1}}$.
Choosing $\rho_0$ small enough, we get the proof.\ $\Box$

Using the method in \cite{JOS}, we have the following lemma:

{\bf Lemma 3.2.}\ \ {\it There exists $v_0\in H^1(\mathbb{R}^N)$ such that $\|v_0\|>\rho_{_{_0}}$ and $I(v_0)<0$.}

{\bf Proof.}\ By $(f_4)$, $\frac{f(t)}{t}$ is decreasing for $t>0$. Consider $\phi\in C^\infty_0(\mathbb{R}^N)$ such that
$0\leq\phi(x)\leq 1$, $\phi(x)=1$ for $|x|\leq 1$, $\phi(x)=0$ for $|x|\geq 2$. We have
$$
f(t\phi(x))\geq f(t)\phi(x),
$$
 for any $x\in \mathbb{R}^N, t>0$. Using $(f_3)$, we get
$$\begin{aligned}
I(t\phi)& = \frac{t^2}{2}\int_{\mathbb{R}^{N}}|\nabla \phi|^{2} +
\frac{1}{2}\int_{\mathbb{R}^{N}}V(x)f^{2}(t\phi)-
\frac{1}{2p}
\int_{\mathbb{R}^{N}}\left(I_\alpha *|f(t\phi)|^p\right)|f(t\phi)|^p\\
& \leq
\frac{t^2}{2}\int_{\mathbb{R}^{N}}|\nabla \phi|^{2} +
\frac{t^2}{2}\int_{\mathbb{R}^{N}}V(x)\phi^2-
\frac{1}{2p}f^{2p}(t)
\int_{\mathbb{R}^{N}}\left(I_\alpha *|\phi|^p\right)|\phi|^p\\
& \leq
\frac{t^2}{2}\left(C_1\|\phi\|^2-C_2\frac{f^4(t)}{t^2}\cdot f^{2p-4}(t)\right).
\end{aligned}$$
By $p>2$ and $(f_6)$, we deduce that $I(t_0\phi)<0$ and $t_0\|\phi\|>\rho_{_{_0}}$ for $t_0$ large enough. Set $v_0=t_0\phi$. Hence $v_0$ is required.\ $\Box$

{\bf Lemma 3.3.}\ \ {\it All Cerami sequences for $I$ at the level $c>0$ are bounded in $H^1(\mathbb{R}^N)$.}

{\bf Proof.} Let $(v_n)\subset H^1(\mathbb{R}^N)$ be a Cerami sequence at the level $c$. Set $w_n:=\frac{f(v_n)}{f'(v_n)}$.
It follows from $(f_4)$ that
$$
\int_{\mathbb{R}^N}|w_n|^2\leq 4\int_{\mathbb{R}^N}|v_n|^2,
$$
$$
\int_{\mathbb{R}^N}|\nabla w_n|^2=\int_{\mathbb{R}^N}\left(1+\frac{2f^2(v_n)}{1+2f^2(v_n)}\right)^2|\nabla v_n|^2
\leq  4\int_{\mathbb{R}^N}|\nabla v_n|^2,
$$
and
$$
|\langle I'(v_n),\ w_n\rangle|\leq C\|I'(v_n)\|\left(\|v_n\|+1\right)\rightarrow 0,\ \ \text{as}\ n\rightarrow \infty.
$$
It follows that $(w_n)\subset H^1(\mathbb{R}^N)$  is bounded.
So
$$\begin{aligned}
c+1& \geq I(v_n)-\frac{1}{2p}\langle I'(v_n), w_n\rangle\\
& =
\frac{1}{2}\int_{\mathbb{R}^{N}}|\nabla v_n|^{2} +
\frac{1}{2}\int_{\mathbb{R}^{N}}V(x)f^{2}(v_n)-
\frac{1}{2p}
\int_{\mathbb{R}^{N}}\left(I_\alpha *|f(v_n^+)|^p\right)|f(v_n^+)|^p\\
& \ \ \ -
\frac{1}{2p}\int_{\mathbb{R}^N}\left(1+\frac{2f^2(v_n)}{1+2f^2(v_n)}\right)|\nabla v_n|^2
-\frac{1}{2p}\int_{\mathbb{R}^{N}}V(x)f^{2}(v_n)\\
& \ \ \ +
\frac{1}{2p}\int_{\mathbb{R}^{N}}\left(I_\alpha *|f(v_n^+)|^p\right)|f(v_n^+)|^p\\
& \geq
\left(\frac{1}{2}-\frac{1}{p}\right)\left(\int_{\mathbb{R}^{N}}|\nabla v_n|^{2} +
\int_{\mathbb{R}^{N}}V(x)f^{2}(v_n)\right).\\
\end{aligned}$$
Since $p>2$, the sequence $\{\int_{\mathbb{R}^{N}}|\nabla v_n|^{2} +
\int_{\mathbb{R}^{N}}V(x)f^{2}(v_n)\}$ is bounded. By the Sobolev embedding theorem and $(f_6)$, we have
$$\begin{aligned}
\int_{\mathbb{R}^{N}}|v_n|^2 & =\int_{\{|v|\leq 1\}}|v_n|^2+\int_{\{|v|>1\}}|v_n|^2\\
& \leq
C_1\int_{\{|v|\leq 1\}}|f(v_n)|^2+\left(\int_{\{|v|>1\}}|v_n|\right)^\theta\left(\int_{\{|v|>1\}}|v_n|^{2^*}\right)^{1-\theta}\\
& \leq
C_1\int_{\mathbb{R}^{N}}f^2(v_n)+\left(\int_{\{|v|>1\}}f^2(v_n)\right)^\theta\left(\int_{\mathbb{R}^{N}}|v_n|^{2^*}\right)^{1-\theta}\\
& \leq
C_2\int_{\mathbb{R}^{N}}V(x)f^2(v_n)+C_3\left(\int_{\mathbb{R}^{N}}V(x)f^2(v_n)\right)^\theta\left(\int_{\mathbb{R}^{N}}|\nabla v_n|^2\right)^{(1-\theta)\frac{2^*}{2}}\\
& \leq
C,
\end{aligned}$$
where $\theta=\frac{2^*-2}{2(2^*-1)}$. Hence $(v_n)$ is bounded in $H^1(\mathbb{R}^N)$.\ $\Box$

{\bf Lemma 3.4.}\ \ {\it Let $\Omega$ be a domain in $\mathbb{R}^N$. Suppose  $\{g_n\},\ \{h_n\}\subset L^1(\Omega)$ and $h\in L^1(\Omega)$. If
 $$0\leq g_n\leq h_n, \ \ g_n(x)\rightarrow 0, \ \ h_n\rightarrow h \ \ a. e.\ \ x\in \Omega$$
and
 $$\lim\limits_{n\rightarrow\infty}\int_{\Omega}h_n=\int_{\Omega}h,$$
 then $\lim\limits_{n\rightarrow\infty} \int_{\Omega}g_n=0$.}\par

{\bf Proof}. By Fatou's lemma,
$$\begin{aligned}
\int_{\Omega}h & =\int_{\Omega}\liminf\limits_{n\rightarrow\infty}(h_n-g_n)\\
& \leq \liminf\limits_{n\rightarrow\infty}\int_{\Omega}(h_n-g_n)\\
& =\liminf\limits_{n\rightarrow\infty}\int_{\Omega}h_n+\liminf\limits_{n\rightarrow\infty}\int_{\Omega}(-g_n)\\
& =\int_{\Omega}h-\limsup\limits_{n\rightarrow\infty}\int_{\Omega}g_n.
\end{aligned}$$
Therefore, $0\leq\liminf\limits_{n\rightarrow\infty}\int_{\Omega}g_n\leq\limsup\limits_{n\rightarrow\infty}\int_{\Omega}g_n\leq 0$, that is, \ $\lim\limits_{n\rightarrow\infty} \int_{\Omega}g_n=0$.\ $\Box$

In the following, we always assume that $\{v_n\}\subset H^1(\mathbb{R}^N)$ is a Cerami sequence for $I$ at the level $c>0$. By the preceding lemma,\ $\{v_n\}$ is bounded. We may assume,
going if necessary to a subsequence, $v_n\rightharpoonup v\in H^1(\mathbb{R}^N)$,\ $v_n(x)\rightarrow v(x)$ a.e. $x\in \mathbb{R}^N$ and $v_n\rightarrow v$ in $L^q_{loc}(\mathbb{R}^N)$ for all $q\in [2, 2^*)$. We have the following Lemma 3.5-3.8.

{\bf Lemma 3.5.}\ \ {\it If $\int_{\mathbb{R}^N}|f(v_n)|^2\rightarrow \int_{\mathbb{R}^N}|f(v)|^2$\ as $n\rightarrow \infty$, then
$\|v_n-v\|\rightarrow 0$.}

{\bf Proof}. The proof of Lemma 3.5 will be carried out in a series of steps.

{\bf Step 1.}\ $\int_{\mathbb{R}^N}V(x)|f(v_n-v)|^2\rightarrow 0$\ as $n\rightarrow \infty$.

By $(f_3)$,\ $\{f(v_n)\}$ is bounded in $L^2(\mathbb{R}^N)$. We can assume $f(v_n)\rightharpoonup f(v)$ in $L^2(\mathbb{R}^N)$, and so
\ $\int_{\mathbb{R}^N}|f(v_n)-f(v)|^2\rightarrow 0$. Therefore, by $(V_1)$ or $(V_1)$ , one has
$$\begin{aligned}
\left|\int_{\mathbb{R}^N}V(x)|f(v_n)|^2-\int_{\mathbb{R}^N}V(x)|f(v)|^2\right| & \leq \int_{\mathbb{R}^N}V(x)\left||f(v_n)|^2
-|f(v)|^2\right|\\
& \leq V_\infty\int_{\mathbb{R}^N}\left||f(v_n)|^2
-|f(v)|^2\right|\rightarrow 0.
\end{aligned}$$
By $(f_{8})$, there is $C>0$ such that $f^2(2t)\leq Cf^2(t)$. Since $f^2(t)$ is convex,
$$\begin{aligned}
V(x)f^2(v_n-v) & \leq
V(x)f^2\left(\frac{1}{2}\left(2v_n-2v\right)\right)\\
& \leq
\frac{1}{2}V(x)f^2(2v_n)+\frac{1}{2}V(x)f^2(-2v)\\
& \leq
\frac{1}{2}V(x)f^2(2v_n)+\frac{1}{2}V(x)f^2(2v)\\
& \leq
C\left(V(x)f^2(v_n)+V(x)f^2(v)\right).
\end{aligned}$$
Using the Lemma 3.4,
$$
\lim\limits_{n\rightarrow \infty}\int_{\mathbb{R}^N}V(x)f^2(v_n-v)=0.
$$

{\bf Step 2.}\ For any $q\in [2, 2^*)$, $\int_{\mathbb{R}^N}|v_n-v|^q\rightarrow 0$ as $n\rightarrow \infty$.

Check the proof of the Lemma 3.3. We have
$$\begin{aligned}
& \ \ \ \ \int_{\mathbb{R}^N}|v_n-v|^2\\
& \leq C_1\int_{\mathbb{R}^{N}}V(x)f^2(v_n-v)
+C_2\left(\int_{\mathbb{R}^{N}}V(x)f^2(v_n-v)\right)^\theta\left(\int_{\mathbb{R}^{N}}|\nabla (v_n-v)|^2\right)^{(1-\theta)\frac{2^*}{2}},
\end{aligned}$$
where $\theta=\frac{2^*-2}{2(2^*-1)}$. Since $(v_n)$ is bounded in $H^1(\mathbb{R}^N)$,\ $\int_{\mathbb{R}^N}|v_n-v|^2\rightarrow 0$ as $n\rightarrow \infty$. It follows from interpolation inequality that $\int_{\mathbb{R}^N}|v_n-v|^q\rightarrow 0$ for any $q\in [2, 2^*)$.

{\bf Step 3.}\ $\|v_n-v\|\rightarrow 0$ as $n\rightarrow 0$.

Using  (2.1),\ $(f_5)$, $(f_{7})$ and the H\"older's inequality, we have
$$\begin{aligned}
&\ \ \ \ \left|\int_{\mathbb{R}^N}\left(I_\alpha*|f(v_n^+)|^p\right)|f(v_n^+)|^{p-2}f(v_n^+)f'(v_n^+)(v_n-v)\right|\\
&\leq C\int_{\mathbb{R}^N}\left(I_\alpha*|v_n|^{p/2}\right)|v_n|^{\frac{p}{2}-1}\left|v_n-v\right|\\
& \leq
C\left(\int_{\mathbb{R}^N}|v_n|^{\frac{p}{2}r}\right)^{1/r}\left(\int_{\mathbb{R}^N}|v_n|^{\frac{p-2}{2}r}|v_n-v|^r\right)^{1/r}\\
& \leq
C\left(\left(\int_{\mathbb{R}^N}|v_n|^{\frac{p-2}{2}r\cdot\frac{p}{p-2}}\right)^{\frac{p-2}{p}}
\left(\int_{\mathbb{R}^N}|v_n-v|^{r\cdot\frac{p}{2}}\right)^{2/p}\right)^{1/r}\\
& \leq
C\left(\int_{\mathbb{R}^N}|v_n-v|^{r\cdot\frac{p}{2}}\right)^{2/(pr)}\rightarrow 0,
\end{aligned}$$
where $\frac{2}{r}-\frac{\alpha}{N}=1$.
Since $\|I'(v_n)\|\rightarrow 0$ and $\{v_n-v\}$ is bounded,
$$
\begin{aligned}
\langle I'(v_n),\ v_n-v\rangle & =\int_{\mathbb{R}^{N}}\nabla v_n\nabla
(v_n-v) + \int_{\mathbb{R}^{N}}V(x)f(v_n)f'(v_n)(v_n-v)\\
& -\int_{\mathbb{R}^{N}}\left(I_\alpha *|f(v_n^+)|^p\right)|f(v_n^+)|^{p-2}f(v_n^+)f'(v_n^+)(v_n-v)
\rightarrow 0.
\end{aligned}
$$
Further,
$$
\left|\int_{\mathbb{R}^{N}}V(x)f(v_n)f'(v_n)(v_n-v)\right|\leq V_\infty \int_{\mathbb{R}^{N}}|v_n||v_n-v|\rightarrow 0.
$$
Hence
$$
\int_{\mathbb{R}^{N}}\nabla v_n\nabla
(v_n-v)\rightarrow 0
$$
and
$$
\int_{\mathbb{R}^{N}}|\nabla(v_n-v)|^2=\int_{\mathbb{R}^{N}}\nabla v_n\nabla
(v_n-v)-\int_{\mathbb{R}^{N}}\nabla v\nabla
(v_n-v)\rightarrow 0.
$$
From arguments above, we get that $\|v_n-v\|\rightarrow 0$ as $n\rightarrow \infty$.\ $\Box$

{\bf Lemma 3.6}\ \ {\it Up to a subsequence,  $A:=\lim\limits_{n\rightarrow\infty}\int_{\mathbb{R}^N}|f(v_n)|^2>0$.}\par

{\bf Proof.}\ We suppose, by contradiction, that $A=0$. By the Lemma 3.5,\ $v_n\rightarrow 0$ in $H^1(\mathbb{R}^N)$. Hence
$$
\int_{\mathbb{R}^N}\left(I_\alpha*|f(v_n^+)|^p\right)|f(v_n^+)|^p\leq C\left(\int_{\mathbb{R}^N}|v_n|^{\frac{pr}{2}}\right)^{\frac{2}{r}}\rightarrow 0
$$
where $\frac{2}{r}-\frac{\alpha}{N}=1$.
Since
$$\begin{aligned}\langle I'(v_n),\ \frac{f(v_n)}{f'(v_n)}\rangle
& =
\int_{\mathbb{R}^N}\left(1+\frac{2f^2(v_n)}{1+2f^2(v_n)}\right)|\nabla v_n|^2
+\int_{\mathbb{R}^{N}}V(x)f^{2}(v_n)\\
& -
\int_{\mathbb{R}^{N}}\left(I_\alpha *|f(v_n^+)|^p\right)|f(v_n^+)|^p\rightarrow 0,
\end{aligned}$$
$$
\int_{\mathbb{R}^N}|\nabla v_n|^2
+\int_{\mathbb{R}^{N}}V(x)f^{2}(v_n)\rightarrow 0.
$$
It follows that
$$
c+o_n(1)=I(v_n)=\frac{1}{2}\left(\int_{\mathbb{R}^N}|\nabla v_n|^2
+V(x)f^{2}(v_n)\right)-\frac{1}{2p}\int_{\mathbb{R}^{N}}\left(I_\alpha *|f(v_n^+)|^p\right)|f(v_n^+)|^p\rightarrow 0,
$$
a contradiction. The proof is completed.\ $\Box$

{\bf Lemma 3.7}\ \ {\it Up to a subsequence, there exist $R,\ \beta>0$ and  $\{x_n\}\subset \mathbb{R}^N$ such that

$$
\liminf\limits_{n\rightarrow +\infty}\int_{B_R(x_n)}|v_n|^2\geq\beta.
$$
}
{\bf Proof.}\ \ By the Lemma 3.6,  up to a subsequence,  one has $A:=\lim\limits_{n\rightarrow\infty}\int_{\mathbb{R}^N}|f(v_n)|^2>0$. If Lemma 3.7 is false, then
 it follows from the Lemma 1.21 in \cite{M1} that, up to a subsequence,
$$
v_n\rightarrow 0\ \ \text{in}\ \ L^2(\mathbb{R}^N).
$$
Hence
$$
0<A=\lim\limits_{n\rightarrow\infty}\int_{\mathbb{R}^N}|f(v_n)|^2\leq \lim\limits_{n\rightarrow\infty}\int_{\mathbb{R}^N}|v_n|^2=0,
$$
a contradiction. This completes the proof.\ $\Box$

{\bf Lemma 3.8}\ \ {\it $\langle I'(v),\ \varphi\rangle=0$ for any $\varphi\in C^\infty_{0}(\mathbb{R}^N)$.}

{\bf Proof.}\ \ For any $\varphi\in C^\infty_{0}(\mathbb{R}^N)$,  the support of $\varphi$ is contained in $B_{R_0}(0)$ for some $R_0>0$.
Hence
$$
\begin{aligned}
& \ \ \ |\langle I'(v_n)-I'(v),\ \varphi\rangle|\\
 & \leq \left|\int_{\mathbb{R}^{N}}\nabla (v_n-v)\nabla
\varphi\right|\\
& + \left|\int_{\mathbb{R}^{N}}V(x)\left(f(v_n)f'(v_n)-f(v)f'(v)\right)\varphi\right|\\
& +\left|\int_{\mathbb{R}^{N}}\left[\left(I_\alpha *|f(v_n^+)|^p\right)|f(v_n^+)|^{p-1}f'(v_n^+)-\left(I_\alpha *|f(v^+)|^p\right)|f(v^+)|^{p-1}f'(v^+)\right]\varphi\right|\\
& :=I_1+I_2+I_3.
\end{aligned}
$$

For $I_1:=\left|\int_{\mathbb{R}^{N}}\nabla (v_n-v)\nabla
\varphi\right|$,\ since $v_n\rightharpoonup v$ in $H^1(\mathbb{R}^N)$,
$I_1\rightarrow 0$ as $n\rightarrow \infty$.\par
For $I_2:= \left|\int_{\mathbb{R}^{N}}V(x)\left(f(v_n)f'(v_n)-f(v)f'(v)\right)\varphi\right|$,\ by $(f_2)$ and $(f_3)$,\ we have $$\left|f(v_n)f'(v_n)-f(v)f'(v)\right|^2\leq 2\left(\left|f(v_n)f'(v_n)\right|^2+\left|f(v)f'(v)\right|^2\right)\leq 2\left|v_n\right|^2+2\left|v\right|^2.$$
By $v_n\rightarrow v$ in $L^2_{loc}(\mathbb{R}^N)$ and the Lemma 3.4, we obtain
$$
\lim\limits_{n\rightarrow \infty}\int_{B_{R_0}(0)}\left|\left(f(v_n)f'(v_n)-f(v)f'(v)\right)\right|^2=0.
$$
Using the H\"older inequality, we have
$$
\begin{aligned}
I_2& \leq V_\infty\int_{B_{R_0}(0)}\left|f(v_n)f'(v_n)-f(v)f'(v)\right||\varphi|\\
& \leq V_\infty\left(\int_{B_{R_0}(0)}\left|f(v_n)f'(v_n)-f(v)f'(v)\right|^2\right)^{\frac{1}{2}}\left(\int_{B_{R_0}(0)}|\varphi|^2\right)^{\frac{1}{2}}\rightarrow 0
\end{aligned}
$$
as $n\rightarrow \infty$.\par
Moreover,
$$\begin{aligned}
I_3 := & \left|\int_{\mathbb{R}^{N}}\left[\left(I_\alpha *|f(v_n^+)|^p\right)|f(v_n^+)|^{p-1}f'(v_n^+)-\left(I_\alpha *|f(v^+)|^p\right)|f(v^+)|^{p-1}f'(v^+)\right]\varphi\right|\\
 \leq &
\int_{\mathbb{R}^{N}}\left(I_\alpha *|f(v_n^+)|^p\right)|\left||f(v_n^+)|^{p-1}f'(v_n^+)-|f(v^+)|^{p-1}f'(v^+)\right||\varphi|\\
& +
\left|\int_{\mathbb{R}^{N}}\left(I_\alpha *|f(v_n^+)|^p\right)|f(v^+)|^{p-1}f'(v^+)\varphi
-\int_{\mathbb{R}^{N}}\left(I_\alpha *|f(v^+)|^p\right)|f(v^+)|^{p-1}f'(v^+)\varphi\right|\\
& :=
J_1+J_2.
\end{aligned}$$
For $r=\frac{2N}{N+\alpha}$, by $(f_5)$ and $(f_{7})$,
$$\begin{aligned}
& \left||f(v_n^+)|^{p-1}f'(v_n^+)-|f(v^+)|^{p-1}f'(v^+)\right|^{r\cdot\frac{p}{p-2}}\\
& \leq
C_1 \left(\left||f(v_n^+)|^{p-2}f(v_n^+)f'(v_n^+)\right|^{r\cdot\frac{p}{p-2}}+\left||f(v^+)|^{p-2}f(v^+)f'(v^+)\right|^{r\cdot\frac{p}{p-2}}\right)\\
& \leq
C_2 \left(\left||f^2(v_n^+)|^{\frac{p-2}{2}}\right|^{r\cdot\frac{p}{p-2}}+\left||f^2(v^+)|^{\frac{p-2}{2}}\right|^{r\cdot\frac{p}{p-2}}\right)\\
& \leq
C_3 \left(|v_n|^{r\cdot\frac{p}{2}}+|v|^{r\cdot\frac{p}{2}}\right)\\
\end{aligned}$$

Since $\frac{2(N+\alpha)}{N}\leq p<\frac{2(N+\alpha)}{N-2}$,\ $\frac{rp}{2}\in [2,\ 2^*)$. By $v_n\rightarrow v$ in $L^{\frac{rp}{2}}_{loc}(\mathbb{R}^N)$ and the Lemma 3.4 again, we obtain
$$
\lim\limits_{n\rightarrow \infty}\int_{B_{R_0}(0)}\left||f(v_n^+)|^{p-1}f'(v_n^+)-|f(v^+)|^{p-1}f'(v^+)\right|^{r\cdot\frac{p}{p-2}}=0.
$$

By the boundedness of $(v_n)$, the H\"older inequality and (2.1), take $n\rightarrow \infty$,
$$\begin{aligned}
J_1& =\int_{\mathbb{R}^{N}}\left(I_\alpha *|f(v_n^+)|^p\right)|\left||f(v_n^+)|^{p-1}f'(v_n^+)-|f(v^+)|^{p-1}f'(v^+)\right||\varphi|\\
& \leq
C_1\left(\int_{\mathbb{R}^{N}}|v_n|^{\frac{p}{2}r}\right)^{\frac{1}{r}}\left(\int_{\mathbb{R}^{N}}|f(v_n^+)|^{p-1}f'(v_n^+)-|f(v^+)|^{p-1}f'(v^+)|^r|\varphi|^r\right)^{\frac{1}{r}}\\
& \leq
C_2 \left(\int_{B_{R_0}(0)}|f(v_n^+)|^{p-1}f'(v_n^+)-|f(v^+)|^{p-1}f'(v^+)|^r|\varphi|^r\right)^{\frac{1}{r}}\\
& \leq
C_3 \left(\int_{B_{R_0}(0)}|f(v_n^+)|^{p-1}f'(v_n^+)-|f(v^+)|^{p-1}f'(v^+)|^{r\cdot\frac{p}{p-2}}\right)^{\frac{p-2}{pr}}\left(\int_{B_{R_0}(0)}|\varphi|^{r\cdot\frac{p}{2}}\right)^{\frac{2}{pr}}\\
& \leq
C_4 \left(\int_{B_{R_0}(0)}|f(v_n^+)|^{p-1}f'(v_n^+)-|f(v^+)|^{p-1}f'(v^+)|^{r\cdot\frac{p}{p-2}}\right)^{\frac{p-2}{pr}}\rightarrow 0
\end{aligned}$$
where $r=\frac{2N}{N+\alpha}$ is given in Remark 2.2-(2).

For $r=\frac{2N}{N+\alpha}$, by $\frac{2(N+\alpha)}{N}\leq p<\frac{2(N+\alpha)}{N-2}$,\ $(f_{7})$ and the H\"older inequality, we have
$$\begin{aligned}
\int_{\mathbb{R}^{N}}||f(v^+)|^{p-1}f'(v^+)\varphi|^r& \leq C\int_{\mathbb{R}^{N}}|f^2(v^+)|^{\frac{p-2}{2}\cdot r}|\varphi|^r\\
& \leq C\int_{\mathbb{R}^{N}}|v|^{\frac{p-2}{2}\cdot r}|\varphi|^r\\
& \leq C\left(\int_{\mathbb{R}^{N}}|v|^{\frac{p-2}{2}\cdot r\cdot\frac{p}{p-2}}\right)^{\frac{p-2}{p}}\left(\int_{\mathbb{R}^{N}}|\varphi|^{r\cdot\frac{p}{2}}\right)^{\frac{2}{p}}\\
& = C\left(\int_{\mathbb{R}^{N}}|v|^{\frac{p}{2}\cdot r}\right)^{\frac{p-2}{p}}\left(\int_{\mathbb{R}^{N}}|\varphi|^{r\cdot\frac{p}{2}}\right)^{\frac{2}{p}}\\
& = C \left|v\right|^{\frac{(p-2)r}{2}}_{\frac{pr}{2}}\left|\varphi\right|^{r}_{\frac{pr}{2}}.
\end{aligned}$$
It follows from $\frac{rp}{2}\in [2,\ 2^*)$ that $|f(v^+)|^{p-1}f'(v^+)\varphi\in L^r(\mathbb{R}^N)$.

In order to prove $J_2\rightarrow 0$, we use an argument which is partly an adaptation of the proof of the Proposition 2.2 in \cite{VJ1}. Set a linear functional $$T(u):=\int_{\mathbb{R}^{N}}\left(I_\alpha *u\right)|f(v^+)|^{p-1}f'(v^+)\varphi.$$ Then, by (2.1),
$T:L^{r}(\mathbb{R}^N)\rightarrow \mathbb{R}$, where $r=\frac{2N}{N+\alpha}$, is a continuous linear functional, that is,
$$\begin{aligned}
|T(u)|& \leq C\left(\int_{\mathbb{R}^{N}}|u|^r\right)^{\frac{1}{r}}\left(\int_{\mathbb{R}^{N}}\left||f(v^+)|^{p-1}f'(v^+)\varphi\right|^r\right)^{\frac{1}{r}}.
\end{aligned}$$
As $(v_n)$ is bounded in $H^1(\mathbb{R}^N)$ and $|f(v_n^+)|^{pr}\leq |v_n|^{\frac{pr}{2}}$, the sequence $(|f(v_n^+)|^p)$ is bounded in $L^r(\mathbb{R}^N)$. We may assume, going if necessary to a subsequence, $|f(v_n^+)|^p\rightharpoonup |f(v^+)|^p$ in $L^r(\mathbb{R}^N)$. Then $T(|f(v_n^+)|^p)\rightarrow T(|f(v^+)|^p)$ as $n\rightarrow \infty$, that is,\ $$J_2=\left|\int_{\mathbb{R}^{N}}\left(I_\alpha *|f(v_n^+)|^p\right)|f(v^+)|^{p-1}f'(v^+)\varphi
-\int_{\mathbb{R}^{N}}\left(I_\alpha *|f(v^+)|^p\right)|f(v^+)|^{p-1}f'(v^+)\varphi\right|\rightarrow 0.$$

So $I_3=J_1+J_2\rightarrow 0$ as $n\rightarrow \infty$. In a summary, up to a subsequence, we prove that $\langle I'(v_n)-I'(v),\ \varphi\rangle\rightarrow 0$ as $n\rightarrow \infty$. Since $\langle I'(v_n),\ \varphi\rangle\rightarrow 0$, we have
$$
\langle I'(v),\ \varphi\rangle=0.\ \Box
$$

{\bf Proof of Theorem 1.1}
As a consequence of the Lemmas 3.1 and 3.2, for the constant
$$
c_{_0}=\inf_{\gamma\in \Gamma}\sup_{t\in[0, 1]}I(\gamma (t))>0,
$$
where
$$
\Gamma=\{\gamma\in C\left([0, 1], H^1(\mathbb{R}^N)\right): \gamma(0)=0, I(\gamma (1))<0\}.
$$
Hence, by  the Theorem 6.3 in \cite{ZFC}, there exists a Cerami sequence $(v_n)$ in $H^1(\mathbb{R}^N)$ at the level $c_{_0}$, that is,
$$
I(v_n)\rightarrow c_{_0}\ \ \text{and}\ \ \left(1+\|v_n\|\right)\|I'(v_n)\|\rightarrow 0,\ \text{as}\ n\rightarrow \infty.
$$
By the Lemma 3.3, the sequence $\{v_n\}$ is bounded.\ Hence, up to a subsequence, one has \ $v_n\rightharpoonup v\in H^1(\mathbb{R}^N)$,\ $v_n(x)\rightarrow v(x)$ a.e. $x\in \mathbb{R}^N$ and $v_n\rightarrow v$ in $L^q_{loc}(\mathbb{R}^N)$ for all $q\in [2, 2^*)$. Hence, by the Lemma 3.8,\ $\langle I'(v),\ \varphi\rangle=0$ for any $\varphi\in C^\infty_0(\mathbb{R}^N)$, that is, $v$ is a weak solution of (1.3). We must prove that $v$ is nontrivial. For this, we follow the idea in \cite{JOS},\ \cite{JYZ2} and \cite{ML} to complete the proof of Theorem 1.1.\par

By the Lemma 3.7, up to a subsequence, there exist $R,\ \beta>0$ and  $\{x_n\}\subset \mathbb{R}^N$ such that
$$
\liminf\limits_{n\rightarrow +\infty}\int_{B_R(x_n)}|v_n|^2\geq\beta.
$$

If $(V_1)$ holds, we may assume that $\{x_n\}$  is bounded.  Then there exists $\rho>0$ such that $B_R(x_n)\subset B_\rho(0)$ for all $n$. Hence
$$
\int_{B_\rho(0)} |v|^2=\liminf\limits_{n\rightarrow +\infty}\int_{B_\rho(0)}|v_n|^2\geq\liminf\limits_{n\rightarrow +\infty}\int_{B_R(x_n)}|v_n|^2\geq\beta>0.
$$
It follows that $v$ is nontrivial. It is easy to see that $v\geq 0$ in $\mathbb{R}^N$. Hence $v\in H^1(\mathbb{R}^N)$ is a nontrivial, nonnegative, weak solution of Eq. $(1.3)$. By the Lemma 2.5,\ $v>0$ in $\mathbb{R}^N$.\par

If $(V_2)$ holds,  we assume, by contradiction, that $v\equiv 0$.  Consider the following two limit functionals:
$$
I_\infty(v)=\frac{1}{2}\int_{\mathbb{R}^{N}}(|\nabla v|^{2} +
V_\infty f^{2}(v))-
\frac{1}{2p}
\int_{\mathbb{R}^{N}}\left(I_\alpha *|f(v^+)|^p\right)|f(v^+)|^p
$$
and
$$
J_\infty(u)=\frac{1}{2}\int_{\mathbb{R}^{N}}[(1+2u^{2})|\nabla u|^{2} +
V_\infty u^{2}] -\frac{1}{2p}
\int_{\mathbb{R}^{N}}\left(I_\alpha *|u^+|^p\right)|u^+|^p,
$$
where $u=f(v)$. Define
$$
c_\infty=\inf_{\gamma\in \Gamma}\sup_{t\in[0, 1]}I_\infty(\gamma (t))>0,
$$
where
$$
\Gamma=\{\gamma\in C\left([0, 1], H^1(\mathbb{R}^N)\right):\ \gamma(0)=0,\ I_\infty(\gamma (1))<0\}.
$$
Notice that  $V_0<V_{\infty}$, we have
$$c_0<c_\infty.  \eqno(3.2)$$

To complete the proof of Theorem 1.1, we divide into the following four lemmas.\par

{\bf Lemma 3.9}\ \ {\it $
I_\infty(v_n)\rightarrow c_0\ \ \text{and}\ \ \left(1+\|v_n\|\right)\|I'_\infty(v_n)\|\rightarrow 0,\ \text{as}\ n\rightarrow \infty.
$
}

{\bf Proof.}\ \ Notice that $(v_n)$ is bounded in $H^1(\mathbb{R}^N)$, there exists $M_1>0$ such that $M_1>2V_\infty$ and $M_1>\int_{\mathbb{R}^N}f^2(v_n)$. Since $v_n\rightarrow v=0$ in $L^q_{loc}(\mathbb{R}^N)$ for all $q\in [2, 2^*)$ and $V(x)\leq V_{\infty}:=\lim\limits_{|y|\rightarrow
 \infty}V(y)<\infty$ for all $x\in \mathbb{R}^{N}$. For every $\epsilon>0$, there is $M>0$ such that, for $n$ large enough, one has
$$
0\leq V_\infty-V(x)<\frac{\epsilon}{2M_1},\ \ \forall \ |x|\geq M,
$$
and
$$
\int_{B_M(0)}|v_n|^2<\frac{\epsilon}{4V_\infty}.
$$
Hence
$$\begin{aligned}
0& \leq \int_{\mathbb{R}^N}V_\infty f^2(v_n)-\int_{\mathbb{R}^N}V(x) f^2(v_n)\\
& =
\int_{\mathbb{R}^N\backslash B_{M}(0)}(V_\infty-V(x)) f^2(v_n)+\int_{B_{M}(0)}(V_\infty-V(x))f^2(v_n)\\
& <
M_1\cdot\frac{\epsilon}{2M_1}+2V_\infty\cdot\frac{\epsilon}{4V_\infty}\\
& =\epsilon,
\end{aligned}$$
and
$$
|I_\infty(v_n)-I(v_n)|=\int_{\mathbb{R}^N}V_\infty f^2(v_n)-\int_{\mathbb{R}^N}V(x) f^2(v_n)\rightarrow 0
$$
as $n\rightarrow +\infty$.
Similarly,
$$\begin{aligned}
\|I'_\infty(v_n)-I'(v_n)\|& =\sup\limits_{\|\phi\|=1}|\langle I'_\infty(v_n)-I'(v_n),\ \phi\rangle|\\
& =
\sup\limits_{\|\phi\|=1}\left|\int_{\mathbb{R}^N}(V_\infty-V(x))f(v_n)f'(v_n)\phi\right|\rightarrow 0
\end{aligned}$$
as $n\rightarrow +\infty$.

It follows that $\{v_n\}$ is also a Cerami sequence of $I_\infty$ at the level $c_0$.\ $\Box$

{\bf Lemma 3.10}\ \ {\it Let $p>2,\ a>0,\ b\geq 0,\ c>0$ and $h(t):=a+bt^2-ct^{2p-2}$ for $t\geq 0$. Then there exists a unique $t_0>0$ such that
$$h(t_0)=0,\ h(t)>0 \ \  \text{for}\ 0\leq t<t_0\ \ \text{and}\ h(t)<0\ \ \text{for}\  t>t_0.$$}
The proof of Lemma 3.10 is standard.\par

{\bf Lemma 3.11.}\ \ {\it Let $v_0\in H^1(\mathbb{R}^N)$,\ $u_0=f(v_0)$ such that
$$a:=\int_{\mathbb{R}^N}|\nabla u_0|^2+\int_{\mathbb{R}^N}V_\infty u_0^2>0,$$
$$b:=4\int_{\mathbb{R}^N}u_0^2|\nabla u_0|^2>0$$
and
$$c:=\int_{\mathbb{R}^N}\left(I_\alpha*|u_0^+|^p\right)|u_0^+|^p>0.$$
Then there exist $t_1>t_0>0$ such that
$$
J_\infty(t_0u_0)>J_\infty(tu_0),\ \ \ \forall \  t\in [0,+\infty)\backslash \{t_0\},
$$
and
$$J_\infty(t_1u_0)<0.$$}
{\bf Proof.}\ \ By the definitions of $I_\infty(v)$ and $J_\infty(u)$, we know that
$$
I_\infty(v_0)=J_\infty(u_0)=\frac{1}{2}\int_{\mathbb{R}^{N}}[(1+2u_0^{2})|\nabla u_0|^{2} +
V_\infty u_0^{2}] -\frac{1}{2p}
\int_{\mathbb{R}^{N}}\left(I_\alpha *|u_0^+|^p\right)|u_0^+|^p.
$$
Set $$g(t)=J_\infty(tu_0)=\frac{1}{2}\int_{\mathbb{R}^{N}}[(1+2t^2u_0^{2})|\nabla (tu_0)|^{2} +
V_\infty t^2u_0^{2}] -\frac{1}{2p}
\int_{\mathbb{R}^{N}}\left(I_\alpha *|tu_0^+|^p\right)|tu_0^+|^p.$$
Then
$$
\frac{d}{dt}(g(t))=t\left(a+bt^2-ct^{2p-2}\right).
$$
By $p>2$ and the Lemma 3.10, there is an unique $t_0>0$ such that $g'(t_0)=0$,\ $g'(t)>0$ for $0<t<t_0$ and $g'(t)<0$ for $t>t_0.$  Hence
$$
J_\infty(t_0u_0)>J_\infty(tu_0),\ \ \ \forall \ t\in [0,+\infty)\backslash \{t_0\}.
$$
Since $p>2$, it is easy to see that $$J_\infty(tu_0)\rightarrow -\infty\ \ \text{as}\ t\rightarrow +\infty.$$
Therefore, there exists $t_1>t_0$ such that $$J_\infty(t_1u_0)<0.\ \Box$$

{\bf Remark 3.12}\ \ Set $$\gamma_0(t)=f^{-1}(tt_1u_0),\ t\in [0,1]. $$
Then
$$
I_\infty(\gamma_0(t))=J_\infty(tt_1u_0)\leq J_\infty(t_0u_0)=J_\infty\left(\frac{t_0}{t_1}\cdot t_1u_0\right)=I_\infty\left(\gamma_0\left(\frac{t_0}{t_1}\right)\right), \ \forall \ t\in [0,1].
$$
It follows that
$$
c_\infty=\inf_{\gamma\in \Gamma}\sup_{t\in[0, 1]}I_\infty(\gamma (t))\leq \sup_{t\in[0, 1]}I_\infty(\gamma_0 (t))=I_\infty\left(\gamma_0\left(\frac{t_0}{t_1}\right)\right)=J_\infty(t_0u_0).
$$

{\bf Lemma 3.13}\ \ {\it $
c_\infty\leq c_0.$}

{\bf Proof.}\ \ By the Lemma 3.6 we have $A:=\lim\limits_{n\rightarrow\infty}\int_{\mathbb{R}^N}|f(v_n)|^2>0$.
Notice that
$$\langle I'_\infty(v_n),\frac{f(v_n)}{f'(v_n)}\rangle=o_n(1).$$
Passing to a subsequence, for $n$ large enough, we get
$$\begin{aligned}
0& <\frac{1}{2}V_\infty A\\
& \leq V_\infty\int_{\mathbb{R}^N}|f(v_n)|^2\\
&= \int_{\mathbb{R}^N}V_\infty|f(v_n)|^2\\
& \leq \langle I'_\infty(v_n),\frac{f(v_n)}{f'(v_n)}\rangle+\int_{\mathbb{R}^{N}}\left(I_\alpha *|f(v_n^+)|^p\right)|f(v_n^+)|^p\\
& =
o_n(1)+\int_{\mathbb{R}^{N}}\left(I_\alpha *|f(v_n^+)|^p\right)|f(v_n^+)|^p
\end{aligned}$$
for $n$ large enough. Hence
$$
\int_{\mathbb{R}^{N}}\left(I_\alpha *|f(v_n^+)|^p\right)|f(v_n^+)|^p\geq \frac{1}{3}V_\infty A>0
$$
for large $n$.
By $\langle I'_\infty(v_n),\ \frac{f(v_n)}{f'(v_n)}\rangle=o_n(1)$ again,
we have
$$
\int_{\mathbb{R}^N}\left(1+\frac{2f^2(v_n)}{1+2f^2(v_n)}\right)|\nabla v_n|^2
+\int_{\mathbb{R}^{N}}V_\infty f^{2}(v_n)\\
 -
\int_{\mathbb{R}^{N}}\left(I_\alpha *|f(v_n^+)|^p\right)|f(v_n^+)|^p=o_n(1).
$$
Set $u_n:=f(v_n)$. Then
$$
\int_{\mathbb{R}^N}(1+4u^2_n)|\nabla u_n|^2
+\int_{\mathbb{R}^{N}}V_\infty u^2_n\\
 -
\int_{\mathbb{R}^{N}}\left(I_\alpha *|u_n^+|^p\right)|u_n^+|^p=o_n(1).
$$
Put
$a_n:=\int_{\mathbb{R}^N}|\nabla u_n|^2+\int_{\mathbb{R}^N}V_\infty u_n^2,$
\ $b_n:=4\int_{\mathbb{R}^N}u_n^2|\nabla u_n|^2$
and
$c_n:=\int_{\mathbb{R}^N}\left(I_\alpha*|u_n^+|^p\right)|u_n^+|^p.$
Then
$$
a_n+b_n-c_n=o_n(1).
$$
Furthermore, since $(u_n)$ is bounded in $H^1(\mathbb{R}^N)$,\ $(a_n),\ (b_n),\ (c_n)$ are all bounded. Hence, passing to a subsequence, we can assume that there are $a,\ b,\ c\in [0, +\infty)$ such that $a_n\rightarrow a,\ b_n\rightarrow b,\ c_n\rightarrow c$ as $n\rightarrow +\infty$ and $a+b-c=0$.
Moreover, for $n$ large enough,  one has
$$
c_n=\int_{\mathbb{R}^N}\left(I_\alpha*|u_n^+|^p\right)|u_n^+|^p=\int_{\mathbb{R}^{N}}\left(I_\alpha *|f(v_n^+)|^p\right)|f(v_n^+)|^p\geq \frac{1}{3}V_\infty A>0
$$
and
$$
a_n=\int_{\mathbb{R}^N}|\nabla u_n|^2+\int_{\mathbb{R}^N}
V_\infty u_n^2\geq \int_{\mathbb{R}^N}
V_\infty f(v_n)^2\geq\frac{1}{2}V_\infty A>0,
$$
 It follows that $a>0,\ c>0.$ It follows from the Lemma 3.10 that there exists a unique sequence $(t_n)\subset (0, +\infty)$ such that $a_n+b_nt_n^2-c_nt_n^{2p-2}=0$. Since $c>0$, $(t_n)$ is bounded. We may assume that there is $t\geq 0$ such that $t_n\rightarrow t.$ Then, $a+bt^2-ct^{2p-2}=0$. Since $a+b-c=0$, by Lemma 3.10 again, we get $t=1$.
By the Lemma 3.11,
$$J_\infty(t_nu_n)>J_\infty(tu_n),\ \forall \ t\in [0, +\infty)\backslash \{t_n\}.$$
Hence $c_\infty\leq J_\infty(t_nu_n)$ by Remark 3.12. Further,

$$\begin{aligned}
J_\infty(t_nu_n)-I_\infty(v_n)& =J_\infty(t_nu_n)-J_\infty(u_n)\\
& =\frac{1}{2}\int_{\mathbb{R}^{N}}[(1+2t_n^2u_n^{2})|\nabla (t_nu_n)|^{2} +
V_\infty t_n^2u_n^{2}]\\
& -\frac{1}{2p}
\int_{\mathbb{R}^{N}}\left(I_\alpha *|t_nu_n^+|^p\right)|t_nu_n^+|^p\\
& -\frac{1}{2}\int_{\mathbb{R}^{N}}[(1+2u_n^{2})|\nabla (u_n)|^{2} +
V_\infty u_n^{2}]\\
& +\frac{1}{2p}
\int_{\mathbb{R}^{N}}\left(I_\alpha *|u_n^+|^p\right)|u_n^+|^p\\
& =
\frac{1}{2}a_n(t_n^2-1)+\frac{1}{4}b_n(t_n^4-1)-\frac{1}{2p}c_n(t_n^{2p}-1)\\
& =
o_n(1).
\end{aligned}$$
Hence
$$
c_\infty\leq J_\infty(t_nu_n)=I_\infty(v_n)+o_n(1).
$$
So,
$$
c_\infty\leq c_0.
$$
\ \ $ \Box$ \par

Contrasting the Lemma 3.13 and (3.2), we get a contradiction. It shows that $v$ is nontrivial. As the case $(V_1)$, we know $v>0$. This completes the proof of Theorem 1.1.\ $\Box$

\end{document}